\documentclass[12pt]{article}
\usepackage{amstext,amsfonts,amsmath,amsthm,graphicx,amssymb,amscd,epsfig}

\newtheorem{Theorem}{Theorem}
\newtheorem{Lemma}{Lemma}
\newtheorem{Proposition}{Proposition}

\begin{document}

\title{ Plane Jacobian conjecture for simple polynomials}
\author{Nguyen Van Chau\thanks{Supported in part by the National Basic Program on Natural Science, Vietnam, and ICTP, trieste, Italy.}\\
{\small Institute of Mathematics, 18 Hoang Quoc Viet, 10307 Hanoi,
Vietnam}\\ {\small E-mail: nvchau@math.ac.vn}}

\date{}

\maketitle

\begin{abstract}
A non-zero constant Jacobian polynomial map $F=(P,Q):\mathbb{C}^2
\longrightarrow \mathbb{C}^2$ has a polynomial inverse if the
component $P$ is a simple polynomial, i.e. if, when $P$ extended
to a morphism $p:X\longrightarrow \mathbb{P}^1$ of a
compactification $X$ of $\mathbb{C}^2$, the restriction of $p$ to
each irreducible component $C$ of the compactification divisor $D
= X-\mathbb{C}^2$ is either degree $0$ or $1$.


\end{abstract}

\noindent{\bf 1.} Let $F=(P,Q):\mathbb{C}^2 \longrightarrow
\mathbb{C}^2$ be a polynomial map, $P,Q\in \mathbb{C}[x,y]$, and
denote $JF:=P_xQ_y-P_yQ_x$ the Jacobian of $F$. The mysterious
Jacobian conjecture (JC) (See \cite{essen-book} and \cite{Bass}),
posed first by Keller in 1939 and still open, asserts that $F$ has
a polynomial inverse if the Jacobian $JF$ is a non-zero constant.
In 1979   by an algebraic approach Razas \cite{Razas} proved this
conjecture for the most simple geometrical case when $P$ is a
rational polynomial, i.e the generic fiber of $P$ is a punctured
sphere, and   all fibres $P=c$, $c\in \mathbb{C}$, are
irreducible. In attempt to understand the geometrical nature of
(JC), this case was also reproved  by Heitmann \cite{Heitmann} and
L\^e and Weber \cite{LeWe2} in some other approaches. In fact, as
observed by Neumann and Norbudy in \cite{NeumannNorbudy}, every
rational polynomial with all irreducible fibres is equivalent to
the coordinate polynomial. Most recent,   L\^e in \cite{Le-hanoi}
and \cite{Le-kyoto} present the following observation, which was
announced in the Hanoi conference, 2006, and the Kyoto conference,
2007.

\begin{Theorem}\label{Le1}{\rm (Theorem 3.2  and Corollary 3.8 in \cite{Le-kyoto})}
A non-zero constant Jacobian polynomial map $F=(P,Q)$ has
 a polynomial inverse if $P$ is  a simple rational polynomial.
 \end{Theorem}
\noindent

Here, following \cite{Neumann-simplePol}, a polynomial map
$P:\mathbb{C}^2\longrightarrow \mathbb{C}$ is {\it simple} if,
when extended $P$ to a morphism $ p:X\longrightarrow \mathbb{P}^1$
of a compactification $X$ of $\mathbb{C}^2$, the restriction of
$p$ to each irreducible component $\ell$ of the compactification
divisor $D = X-\mathbb{C}^2$ is either of  degree $0$ or $1$. In
fact, as in the proof in \cite{Le-kyoto} of Theorem 1, if a
component of non-zero constant Jacobian map $F=(P,Q)$ is a simple
rational polynomial, then this component determines a locally
trivial fibration.

In this short paper we would like to present an other explanation
for Theorem 1 from view point of the geometry of the non-proper
valuer set of the map $F$. In fact, we shall prove

\begin{Theorem}\label{main} A non-zero constant Jacobian polynomial map $F=(P,Q)$ has
 a polynomial inverse if $P$ is  a simple  polynomial.
 \end{Theorem}

In any meaning, the  addition condition on the simple polynomial
component in a non-zero constant Jacobian polynomial map may be
viewed as  a kind of ``good" local conditions at infinity, but it
seems to be not a global one. A completed classification of all
simple rational polynomials was presented in
\cite{Neumann-simplePol}.
\bigskip

\noindent{\bf 2.} Given a polynomial map $F=(P,Q)$ of
$\mathbb{C}^2$. Following \cite{Jelonek}, the non-proper value set
$A_F$ of $F$ is the set of all values $a\in\mathbb{C}^2$ such that
there exists a sequence $\mathbb{C}^2 \ni b_i \rightarrow \infty$
with $F(b_i)\rightarrow a$. This set $A_F$ is either empty or an
algebraic curve in $\mathbb{C}^2$ for which every irreducible
component is the image of a non-constant polynomial map from
$\mathbb{C}$ into $\mathbb{C}^2$. Our argument in the proof of
Theorem 2 here is based on  the following facts, that was
presented in \cite{Chau99} and can be reduced from \cite{Cassou}
(see also \cite{Chau2004} and \cite{Chau2007} for other refine
versions).

\begin{Theorem} Support $F=(P,Q)$ is a polynomial map with non-zero constant Jacobian.
If $A_F\neq \emptyset$, then every irreducible components of $A_f$
can be parameterized by polynomial maps $\xi\mapsto
(\varphi(\xi),\psi(\xi))$ with
$$\deg \varphi/\deg \psi=\deg P/\deg Q.\eqno (1) $$
 \end{Theorem}
 \noindent This theorem together with the
Abhyankar-Moh Theorem \cite{shreeram} on the embeddings of the
line to the plane allows us to obtain:

\begin{Theorem}\label{Theo3} A polynomial map $F$ of $\mathbb{C}^2$ must
have singularities if its non-proper value set $A_F$ has an
irreducible component isomorphic to the line.
\end{Theorem}

\noindent A simple proof of Theorem 4  recently presented in
\cite{Chau2007} gives a description in terms of Newton-Puiseux
data how the singularity occurs in this situation.

\bigskip

\noindent{\bf 3.} To use Theorem \ref{Theo3} in the situation of
simple polynomials, at first, we need to describe the the
non-proper value curve $A_F$ in terms of the regular extension of
$F$ in a convenience compatification $X\supset \mathbb{C}^2$.
Given a polynomial $F=(P,Q)$, extend it to a map $F :
\mathbb{P}^2\longrightarrow \mathbb{P}^1\times \mathbb{P}^1$ and
resolve the points of indeterminacy to get a regular map $ f=(p,q)
: X \longrightarrow\mathbb{P}^1\times\mathbb{P}^1$that coincides
with $F=(P,Q)$ on $\mathbb{C}^2\subset X$. We call $D = X -
\mathbb{C}^2$ the divisor at infinity. The divisor $D$ consists of
a connected union of rational curves isomorphism to $\mathbb{P}^1$
and the dual graph of the divisor $D$ is a tree. An irreducible
component $\ell$ of $D$ is a {\it horizontal }curve of $P$ (or
$Q$) if the restriction of $p$ (res. $q$) to $E$ is not a constant
mapping.  An irreducible component $\ell$ of $D$ is a {\it
dicritical} curve of $F$ if the restriction of $f$ to $\ell$ is
not a constant mapping. A dicritical curve of $F$ must be a
horizontal curve of $P$ or $Q$. Although the compactification
defined above is not unique, the horizontal curves of $P$ or $Q$
as well as the dicritical curves of $F$ are essentially
independent of choice. Further, by blow-down components of
self-intersection $-1$ corresponding to linear vertexes or
endpoint in the dual graph of $D$, but not of dicritical curves of
$F$, horizontal curves of $P$ or $Q$, we can work with minimal
compactification $X\supset\mathbb{C}^2$ on which $F$ can be
extended to a regular morphism $f=(p,q):X\longrightarrow
\mathbb{P}^1\times\mathbb{P}^1$.

Let us denote $D_\infty:=f^{-1}((\{\infty\}\times \mathbb{P}^1)
\cup (\mathbb{P}^1 \times \{\infty\}))$. The following description
of the dual graph of the divisor $D$ is well-known (see, for
example, in \cite{Vitushkin}, \cite{Orevkov} and \cite{LeWe1}).

\begin{Proposition}\label{Pro1}

i) The dual graph of the divisor $D$ is a tree;

ii) The dual graph of the curve $D_\infty$ is a tree;

iii) The dual graph of each connected component  of the closure of
$(D-D_\infty)$ is  a linear path of the form
$$
\odot\rightarrow\circ\rightarrow\circ\rightarrow \dots
\circ\rightarrow\circ
$$
in which the beginning vertex $\odot$ is a dicritical curve of $F$
and the next possible vertexes $\circ$  are curves to which the
restriction of $f$ are finite constant mappings.

\end{Proposition}

The following provides a description of the non-proper value set
$A_F$ of $F$ in terms of regular extension of $F$ in a minimal
compatification $X\supset \mathbb{C}^2$.

\begin{Proposition}\label{Pro2}

i)$$A_F=\bigcup_{\mbox {dicritical curves }\ell \mbox{ of
}F}f(\ell )\cap \mathbb{C}^2.\eqno (2)$$

ii) Let $\ell$ be a dicritical curve of $F$. Then, $\ell$ and the
curve $D_\infty$ have an unique common point. Let $\ell^*:=\ell
-D_\infty$. Then, the curve $\ell^*$ is isomorphic to $\mathbb{C}$
and
$$f(\ell^*)=f(\ell )\cap \mathbb{C}^2\eqno (3)$$

iii) $$A_F=\bigcup_{\mbox {dicritical curves }\ell \mbox{ of
}F}f(\ell^* ).\eqno (4)$$

\end{Proposition}

\begin{proof} Conclusion (i) can be easy verified by the definition of
$A_F$ and a simple topological argument. Conclusion (ii) follows
from the fact in Proposition 1 that the dual graph of the divisor
$D$ is a tree. Conclusion (iii) results from (i) and (ii).
\end{proof}

\bigskip
\noindent{\bf 4.} Now, we consider the situation when the
restriction of $p$ to a dicritical curve $\ell$ of $F$ is of
degree $1$.

\begin{Lemma}\label{Lem1}
Let $\ell$ be a dicritical component of $F$. If the restriction of
$p$ to $\ell$ is of  degree $1$, then the image $f(\ell^*)$ is
isomorphic to the line $\mathbb{C}$.
\end{Lemma}

\begin{proof}
Suppose $\ell$ is a dicritical component of $F$ and the degree of
the restriction $p_{|\ell}$ equals $1$. Then, $p_{|\ell}: \ell
\longrightarrow \mathbb{P}^1$ is injective, and hence, is
bijective, since $\ell$ is isomorphic to $\mathbb{P}^1$. This
ensures that the curve $f(\ell^*)$ intersects each line $\{
(u,v)\in \mathbb{C}^2: u=c\}$, $c\in \mathbb{C}$, at an unique
point. Then, the polynomial $H(u,v)$ defining the curve $f(\ell^*)
\subset \mathbb{C}^2$ can be chosen of the form $v+h(u)$,
$h\in\mathbb{C}[u]$. So, the automorphism $A(u,v):=(u, v-h(u))$
maps isomorphically the curve $f(\ell^*)$ onto the line $v=0$.
\end{proof}

\medskip
\noindent{\it Proof of Theorem \ref{main}. } Suppose $F=(P,Q)$
with $JF\equiv c \neq 0$ and $P$ is a simple polynomial. Note that
each dicritical curve of $F$ must be a horizontal curve of $P$ or
$Q$. Since $JF\equiv c \neq 0$ and $P$ is simple, in view of
Theorem \ref{Theo3} and Lemma \ref{Lem1}, a horizontal curve of
$P$ cannot be   a dicritical curve of $F$. So, if $\ell$ is a
dicritical curve of $F$, then $\ell$ must be a horizontal curve of
$Q$ and the restriction $p_{|\ell}$ maps $\ell$ to a finite
constant. Thus, for such $\ell$ the image $f(\ell^*)$ is a line
$u=const.$. The last is impossible again by Theorem \ref{Theo3} as
$JF\equiv c\neq 0$. Hence, $F$ has not any dicritical component.
Then,  $A_F=\emptyset$ by Proposition \ref{Pro2} and $F$ is a
proper map by the definition of $A_F$. Therefore, by simple
connectedness of $\mathbb{C}^2$ the local diffeomorphic map $F$
must be bijective. Thus, $F$ is an automorphism of $\mathbb{C}^2$.

\bigskip

\noindent{\it Acknowledgments}:  The author wishes to thank Prof.
L\^e D\~ung Tr\'ang for his helps and useful discussions on  the
Jacobian problem.

\end{document}